\documentstyle[12pt]{article}
\begin{document}
\title{  An extremal problem on potentially $K_{p,1,1}$-graphic sequences
\thanks{  Project Supported by NNSF of China(10271105), NSF of Fujian,
Science and Technology Project of Fujian, Fujian Provincial
Training Foundation for "Bai-Quan-Wan Talents Engineering" ,
Project of Fujian Education Department and Project of Zhangzhou
Teachers College.}}
\author{{Chunhui Lai}\\
{\small Department of Mathematics}\\{\small Zhangzhou Teachers
College, Zhangzhou} \\{\small Fujian 363000,
 P. R. of CHINA.}\\{\small e-mail: zjlaichu@public.zzptt.fj.cn}}
\date{}
\maketitle
\begin{center}
\begin{minipage}{120mm}
\vskip 0.1in
\begin{center}{\bf Abstract}\end{center}
 {A sequence $S$ is potentially $K_{p,1,1}$ graphical if it has
a realization containing a $K_{p,1,1}$ as a subgraph, where
$K_{p,1,1}$ is a complete 3-partite graph  with partition sizes
$p,1,1$. Let $\sigma(K_{p,1,1}, n)$ denote the smallest degree sum
such that every $n$-term graphical sequence $S$ with
$\sigma(S)\geq \sigma(K_{p,1,1}, n)$ is potentially $K_{p,1,1}$
graphical.  In this paper, we prove that $\sigma (K_{p,1,1},
n)\geq 2[((p+1)(n-1)+2)/2]$ for $n \geq p+2.$  We conjecture that
equality holds for $n \geq 2p+4.$ We prove that this conjecture is
true for $p=3$.}\par
\par
\end{minipage}
\end{center}
 \par
 \section{Introduction}
\par

  If $S=(d_1,d_2,...,d_n)$ is a sequence of
non-negative integers, then it is called  graphical if there is a
simple graph $G$ of order $n$, whose degree sequence ($d(v_1 ),$
$d(v_2 ),$ $...,$ $d(v_n )$) is precisely $S$. If $G$ is such a
graph then $G$ is said to realize $S$ or be a realization of $S$.
A graphical sequence $S$ is potentially $H$ graphical if there is
a realization of $S$ containing $H$ as a subgraph, while $S$ is
forcibly $H$ graphical if every realization of $S$ contains $H$ as
a subgraph. Let $\sigma(S)=d(v_1 )+d(v_2 )+... +d(v_n ),$ and
$[x]$ denote the largest integer less than or equal to $x$.  We
denote  $G+H$ as the graph with $V(G+H)=V(G)\bigcup V(H)$ and
$E(G+H)=E(G)\bigcup E(H)\bigcup \{xy: x\in V(G) , y \in V(H) \}.$
Let $K_k$, and $C_k$ denote a complete graph on $k$ vertices, and
a cycle on $k$ vertices, respectively. Let $K_{p,1,1}$ denote a
complete 3-partite graph  with partition sizes $p,1,1.$\par

Given a graph $H$, what is the maximum number of edges of a graph
with $n$ vertices not containing $H$ as a subgraph? This number is
denoted $ex(n,H)$, and is known as the Tur\'{a}n number. This
problem was proposed for $H = C_4$ by Erd\"os [3] in 1938 and in
general by Tur\'{a}n [12]. In terms of graphic sequences, the
number $2ex(n,H)+2$ is the minimum even integer $l$ such that
every $n$-term graphical sequence $S$ with $\sigma (S)\geq l $ is
forcibly $H$ graphical. Here we consider the following variant:
determine the minimum even integer $l$ such that every $n$-term
graphical sequence $S$ with $\sigma(S)\ge l$ is potentially $H$
graphical. We denote this minimum $l$ by $\sigma(H, n)$. Erd\"os,\
Jacobson and Lehel [4] showed that $\sigma(K_k, n)\ge
(k-2)(2n-k+1)+2$ and conjectured that equality holds. They proved
that if $S$ does not contain zero terms, this conjecture is true
for $k=3,\ n\ge 6$. The conjecture is confirmed in [5],[7],[8],[9]
and [10].
 \par
 Gould,\ Jacobson and
Lehel [5] also proved that  $\sigma(pK_2, n)=(p-1)(2n-2)+2$ for
$p\ge 2$; $\sigma(C_4, n)=2[{{3n-1}\over 2}]$ for $n\ge 4$. Luo
[11] characterized the potentially $C_{k}$ graphic sequence for
$k=3,4,5.$  Yin and Li [13] gave sufficient conditions for a
graphic sequence being potentially $K_{r,s}$-graphic, and
determined $\sigma(K_{r,r},n)$ for $r=3,4.$ Lai [6] proved that
$\sigma (K_4-e, n)=2[{{3n-1}\over 2}]$ for $n\ge 7$.\  In this
paper, we prove that $\sigma (K_{p,1,1}, n)\geq
2[((p+1)(n-1)+2)/2]$ for $n \geq p+2.$  We conjecture that
equality holds for $n \geq 2p+4.$ We prove that this conjecture is
true for $p=3$.\par

\section{ Main results.} \par
{\bf  Theorem 1.} $\sigma (K_{p,1,1}, n)\geq 2[((p+1)(n-1)+2)/2],$
for $n \geq p+2.$
\par
{\bf Proof.}  If $p=1, $ by Erd\"os,\ Jacobson and Lehel [4],
$\sigma (K_{1,1,1}, n)\geq 2n, $ Theorem 1 is true. If $p=2, $ by
Gould,\ Jacobson and Lehel [5], $\sigma (K_{2,1,1}, n)= \sigma
(K_{4}-e, n)\geq \sigma (C_{4}, n)=2[(3n-1)/2], $ Theorem 1 is
true. Then we can suppose that $p\geq 3.$
\par
We first consider odd $p$. If $n$ is odd, let $n=2m+1, $ by
Theorem 9.7 of [2], $K_{2m}$ is the union of one 1-factor $M$ and
$m-1$ spanning cycles $C_{1}^{1}, C_{2}^{1}, ..., C_{m-1}^{1}.$
Let
$$H=C_{1}^{1}\bigcup C_{2}^{1}\bigcup...\bigcup
C_{\frac{p-1}{2}}^{1}+K_{1}$$ Then $H$ is a realization of
$((n-1)^{1}, p^{n-1}),$ where the symbol $x^{y}$ stands for y
consecutive terms $x$.  Since  $K_{p,1,1}$ contains two vertices
of degree $p+1$ while $((n-1)^{1}, p^{n-1})$ only contains one
integer $n-1$ greater than degree $p$,  $((n-1)^{1}, p^{n-1})$ is
not potentially $K_{p,1,1}$ graphic. Thus
$$\sigma (K_{p,1,1}, n)\geq (n-1) + p(n-1)+ 2
= 2[((p+1)(n-1)+2)/2].$$ Next, if $n$ is even, let $n=2m+2, $ by
Theorem 9.6 of [2], $K_{2m+1}$ is the union of  $m$ spanning
cycles $C_{1}^{1}, C_{2}^{1}, ..., C_{m}^{1}.$ Let
$$H=C_{1}^{1}\bigcup C_{2}^{1}\bigcup...\bigcup
C_{\frac{p-1}{2}}^{1}+K_{1}$$ Then $H$ is a realization of
$((n-1)^{1}, p^{n-1}),$ and  we are done as before. This completes
the discussion for odd $p$.
\par
Now we  consider even $p$. If $n$ is odd, let $n=2m+1, $ by
Theorem 9.7 of [2], $K_{2m}$ is the union of one 1-factor $M$ and
$m-1$ spanning cycles $C_{1}^{1}, C_{2}^{1}, ..., C_{m-1}^{1}.$
Let
$$H=M\bigcup C_{1}^{1}\bigcup C_{2}^{1}\bigcup...\bigcup
C_{\frac{p-2}{2}}^{1}+K_{1}$$ Then $H$ is a realization of
$((n-1)^{1}, p^{n-1}),$ and  we are done as before. Next, if $n$
is even, let $n=2m+2, $ by Theorem 9.6 of [2], $K_{2m+1}$ is the
union of  $m$ spanning cycles $C_{1}^{1}, C_{2}^{1}, ...,
C_{m}^{1}.$ Let
$$C_{1}^{1}= x_{1}x_{2}...x_{2m+1}x_{1}$$
$$H=(C_{1}^{1}\bigcup C_{2}^{1}\bigcup...\bigcup
C_{\frac{p}{2}}^{1}+K_{1})-\{x_{1}x_{2}, x_{3}x_{4}, ...,
x_{2m-1}x_{2m}, x_{2m+1}x_{1}\}$$ Then $H$ is a realization of
$((n-1)^{1}, p^{n-2}, (p-1)^{1}).$ It is easy to see that
$((n-1)^{1}, p^{n-2}, (p-1)^{1})$ is not potentially $K_{p,1,1}$
graphic. Thus
$$\sigma (K_{p,1,1}, n)\geq (n-1) + p(n-2)+ p-1 + 2
= 2[((p+1)(n-1)+2)/2].$$
 This completes
the discussion for even $p$, and so finishes the proof of Theorem
1.
\par
{\bf  Theorem 2.} For $n=5$ and $n\geq7$,
    $$\sigma(K_{3,1,1} ,n)=4n-2.$$
For $n=6$, if S is a 6-term graphical sequence with $\sigma(S)
\geq 22$, then either there is a realization of S containing
$K_{3,1,1}$ or $S=(4^{6})$. (Thus $\sigma(K_{3,1,1} ,6)=26$.)
\par
{\bf Proof.} By theorem 1, for $n\geq5,$
$\sigma(K_{3,1,1,}n)\geq2{[((3+1)(n-1)+2)/2]}=4n-2.$ We need to
show that if S is an n-term graphical sequence with
$\sigma(S)\geq4n-2$, then there is a realization of S containing a
$K_{3,1,1}$ (unless $S=(4^{6})$).
 Let $d_{1} \geq d_{2} \geq \cdots \geq d_{n}$, and let $G$ is a realization of S.
 \par
   Case: $n=5$, if a graph has size
$q \geq9$, then clearly it contains a $K_{3,1,1}$, so that
$\sigma(K_{3,1,1} ,5)\leq4n-2$.
\par
  Case: $n=6$, If $\sigma(S)=22$, we first consider $d_{6} \leq2$. Let
  $S^{'}$ be the degree sequence of $G-v_{6}$, so
  $\sigma(S^{'})\geq22-2\times2=18$. Then $S^{'}$ has a realization
  containing a $K_{3,1,1}$. Therefore S has a realization
   containing a $K_{3,1,1}$. Now we consider $d_{6}\geq3$. It is
  easy to see that $S$ is one of $(5^{2},3^{4})$ ,  $(5^{1},4^{2},3^{3})$
  or $(4^{4},3^{2})$. Obviously, all of them are potentially
  $K_{3,1,1}$ -graphic.
   Next, if $\sigma(S)=24$, we first consider $d_{6}\leq3$. Let
  $ S^{'}$ be the degree sequence of $G-v_{6}$, so
   $\sigma(S^{'})\geq24-3\times2=18$. Then $S^{'}$ has a realization
   containing a $K_{3,1,1}$. Therefore S has a realization
   containing a $K_{3,1,1}$. Now we consider $d_{6}\geq4$. It is
   easy to see that $S=(4^{6})$. Obviously, $(4^{6})$ is graphical
   and $(4^{6})$ is not potentially $K_{3,1,1}$ graphic.
   Finally, suppose that $\sigma(S)\geq26$. We first consider
   $d_{6}\leq4$. Let
 $ S^{'}$ be the degree sequence of $G-v_{6}$, so
  $\sigma(S^{'})\geq26-2\times4=18$. Then $S^{'}$ has a realization
  containing a $K_{3,1,1}$. Therefore S has a realization
   containing a $K_{3,1,1}$. Now we consider $d_{6}\geq5$. It is
   easy to see that $S=(5^{6})$. Obviously, $(5^{6})$ is potentially $K_{3,1,1}$-graphic.

  \par
   Case: $n=7$. First we assume that $\sigma(S)=26$. Suppose $d_{7}\leq2$ and let $S^{'}$ be the degree sequence of
    $G-v_{7}$, so $\sigma(S^{'})\geq26-2\times2=22$. Then $S^{'}$ has a realization
  containing a $K_{3,1,1}$ or $S^{'}=(4^{6})$. Therefore S has a realization
   containing a $K_{3,1,1}$ or $S=(5^{1},4^{5},1^{1})$. Obviously, $(5^{1},4^{5},1^{1})$
   is potentially $K_{3,1,1}$ -graphic. In either event, S has a realization
   containing a $K_{3,1,1}.$  Now we
   assume that $d_{7}\geq3$. It is easy to see that S is one of
   $(6^{1},5^{1},3^{5})$,
   $(6^{1},4^{2},3^{4})$, $(5^{2},4^{1},3^{4})$,$(5^{1},4^{3},3^{3})$
   or $(4^{5},3^{2})$. Obviously, all of them are  potentially $K_{3,1,1}$-graphic.
Next,  if $\sigma(S)=28$, Suppose $d_{7}\leq3$. Let $S^{'}$ be the
degree sequence of $G-v_{7}$, so
$\sigma(S^{'})\geq28-3\times2=22$. Then $S^{'}$ has a realization
  containing a $K_{3,1,1}$ or $S^{'}=(4^{6})$. Therefore S has a realization
   containing a $K_{3,1,1}$ or $S=(5^{2},4^{4},2^{1})$. Obviously, $(5^{2},4^{2},2^{1})$
   is potentially $K_{3,1,1}$ -graphic. In either event, S has a realization
   containing a $K_{3,1,1}.$  Now we
   assume that $d_{7}\geq4$, then $S=(4^{7})$. Clearly, $(4^{7})$ has a realization
  containing a $K_{3,1,1}$.
  Finally, suppose that $\sigma(S)\geq30$. If $d_{7}\leq4$. Let $S^{'}$ be the degree sequence
  of $G-v_{7}$, so $\sigma(S^{'})\geq30-2\times4=22$. Then $S^{'}$ has a realization
  containing a $K_{3,1,1}$ or $S^{'}=(4^{6})$. Therefore S has a realization
   containing a $K_{3,1,1}$ or $S=(5^{3},4^{3},3^{1})$.
   Clearly, $(5^{3},4^{3},3^{1})$ has a realization containing a
   $K_{3,1,1}$. In either event, S has a realization
   containing a $K_{3,1,1}.$   Now we consider $d_{7}\geq5$. It is easy to see that
   $\sigma(S)\geq5\times7=35$. Obviously
   $\sigma(S)\geq36$. Clearly, S has a realization containing a $K_{3,1,1}.$
   \par
   We proceed by induction on n. Take $n \geq 8$ and make the
   inductive assumption that for $7 \leq t < n$, whenever $S_{1}$ is
   a t-term graphical sequence such that
          $$\sigma(S_{1})\geq4t-2$$
   then $S_{1}$ has a realization containing a $K_{3,1,1}.$
      Let S be an n-term graphical sequence with $\sigma(S)\geq4n-2$.
      If $d_{n}\leq2$, let $S^{'}$ be the degree sequence of $G-v_{n}$. Then
   $\sigma(S^{'})\geq4n-2-2\times2=4(n-1)-2$. By induction, $S^{'}$ has a realization
   containing a $K_{3,1,1}$. Therefore S has a realization containing a
   $K_{3,1,1}$. Hence, we may assume that $d_{n}\geq3$. By
   Proposition 2 and Theorem 4 of [5] (or Theorem  3.3 of [7] ) S has a
   realization containing a $K_{4}$. By Lemma 1 of [5] ,there is a
   realization G of S with $v_{1},v_{2},v_{3},v_{4}$,  the four
   vertices of highest degree containing a $K_{4}$.
     If $d(v_{2})=3$, then $4n-2 \leq \sigma(S) \leq n-1+3(n-1)=4n-4$. This is a
   contradiction. Hence, we may assume that $d(v_{2}) \geq 4$.
      Let $v_{1}$ be adjacent to $v_{2},v_{3},v_{4},y_{1}$. If $y_{1}$ is adjacent
   to one of $v_{2},v_{3},v_{4}$, then G contains a $K_{3,1,1}$. Hence, we may assume
   that $y_{1}$ is not adjacent to $v_{2},v_{3},v_{4}$.
      Let $v_{2}$ be adjacent to $v_{1},v_{3},v_{4},y_{2}$. If $y_{2}$ is adjacent
   to one of $v_{1},v_{3},v_{4}$, then G contains a $K_{3,1,1}$. Hence, we may assume
   that $y_{2}$ is not adjacent to $v_{1},v_{3},v_{4}$.
       Since $d(y_{1}) \geq d_{n} \geq 3$, there is a new vertex $y_{3}$, such that
   $y_{1}y_{3} \in E(G)$.
   \par
       Case 1: Suppose $y_{3}v_{3} \in E(G)$.
       If $y_{3}v_{4} \in E(G)$, then G contains a $K_{3,1,1}$. Hence, we may assume
   that $y_{3}v_{4}  \notin  E(G)$. Then the edge interchange that removes the edges
   $y_{1}y_{3},v_{3}v_{4}$ and $v_{2}y_{2}$ and inserts the edges
   $y_{1}v_{2},y_{3}v_{4}$ and $y_{2}v_{3}$ produces a realization
   $G^{'}$ of S containing a $K_{3,1,1}$.
       \par
       Case 2: Suppose $y_{3}v_{3} \notin  E(G)$. Then the edge interchange that removes
   the edges $y_{1}y_{3},v_{3}v_{4}$  and $v_{2}y_{2}$ and inserts
   the edges $y_{1}v_{2},y_{3}v_{3}$  and $y_{2}v_{4}$ produces a
   realization $G^{'}$ of S containing a $K_{3,1,1}$.
       This finishes the inductive step, and thus  Theorem 2 is established.
      \par
       We make the following conjecture:
       \par
     {\bf Conjecture.}  $$\sigma(K_{p,1,1} ,n)=2[((p+1)(n-1)+2)/2]$$
   for $n \geq 2p+4$.
   \par
       This conjecture is true for $p=1$, by Theorem 3.5  of [4],
    for $p=2$, by Theorem 1 of [6], and  for $p=3$, by the above Theorem 2.
       \par

 \section*{Acknowledgment}

 The author thanks Prof. Therese Biedl for her valuable suggestions.
 The author thanks the referees for many helpful comments.
 \par

\end{document}